\documentclass[]{interact}

\usepackage{epstopdf}
\usepackage[caption=false]{subfig}
\usepackage{dsfont,bm}
\usepackage[numbers,sort&compress]{natbib}
\bibpunct[, ]{[}{]}{,}{n}{,}{,}
\makeatletter
\def\NAT@def@citea{\def\@citea{\NAT@separator}}
\makeatother

\theoremstyle{plain}

\theoremstyle{definition}

\theoremstyle{remark}

\begin{document}


\title{Multivariate two-sample test statistics based on data depth}

\author{
\name{Yiting Chen\textsuperscript{a} \textsuperscript{b} \thanks{CONTACT Xiaoping Shi. Email: xiaoping.shi@ubc.ca}, Wei Lin \textsuperscript{b}, Xiaoping Shi \textsuperscript{a}}
\affil{\textsuperscript{a}Department of Computer Science, Mathematics, Physics and Statistics, University of British Columbia, Kelowna, Canada V1V 1V7; \\
\textsuperscript{b}Department of Statistics and Actuarial Science, Simon Fraser University, Burnaby, BC, Canada V5A 1S6 }}

\maketitle

\begin{abstract}
Data depth has been applied as a nonparametric measurement for ranking multivariate samples. In this paper, we focus on homogeneity tests to assess whether two multivariate samples are from the same distribution. There are many data depth-based tests for this problem, but they may not be very powerful, or have unknown asymptotic distributions, or have slow convergence rates to asymptotic distributions. Given the recent development of data depth as an important measure in quality assurance, we propose three new test statistics for multivariate two-sample homogeneity tests. The proposed minimum test statistics have simple asymptotic half-normal distribution. We also discuss the generalization of the proposed tests to multiple samples. The simulation study demonstrates the superior performance of the proposed tests. The test procedure is illustrated by two real data examples.
\end{abstract}


\section{Introduction}
Multivariate statistical analysis has been widely applied in many fields in recent years. We focus on the tests for homogeneity multivariate two samples, i.e., $H_0: F=G$ vs $H_\alpha: F \neq G$ for two multivariate samples from $F$ and $G$ distributions, respectively. There are two types of tests: parametric and non-parametric. The parametric tests include multivariate analysis of variance (MANOVA), for the samples under the normality assumption. Typical non-parametric tests include Cram\'er test, Energy distance test, and Depth-based Tests. This paper focuses on non-parametric statistical tests based on data depth. 
The depth function $D(x;F)$ measures the centrality of one point $x$ in distribution $F(x)$ and transforms from space $R^d$ into [0,1] in a d-dimensional space. The existing two-sample depth based tests, such as Depth-based Rank Test \cite{Small2011}, Weighted and Maximum Statistics \cite{Shi2023}, are not very powerful with a relatively low power with small sample size. Hence, we proposed three new test statistics, which are superior than the existing test statistics. We used three depth functions: Mahalanobis depth \cite{liu1993quality}, Spatial depth \cite{Brown58,Gower74}, and Projection depth \cite{Liu92}, which can be implemented via the \textit{R} package \textit{ddalpha}.

Statistical depth has the advantage that it does not need the condition on normality and provides the rank of distributions.  \cite{Serfling2000,liu1993quality} listed the properties of data depth functions: 
\begin{enumerate}
\item Affine invariance: the coordinate system or scales will not influence the depth.
\item Maximality at center: the center point of a distribution has a maximum value, i.e., $D(\mu,F)=\sup_{x\in \mathcal{R}^d} D(x,F)$ and $\mu$ is the center point.
\item Monotonicity relative to deepest point: With center$\mu$ in distribution $F$, as any point x moves farther away from the center, the depth value decreases monotonically, i.e., $D(x,F) <D(\mu+ \alpha(x-\mu),F)$ for any $0<\alpha<1$.
\item Vanishing at infinity: The depth value will go to zero as the $||x||$ goes to infinity.
\end{enumerate}

Under the statistical depth definition, the Q Statistics, quality index, is proposed in \cite{liu1993quality}. Our new proposed test statistics are based on this quality index and used for two-sample problems. It measures the relative ``outlyingness" of one distribution to another. In particular, in $Q(F,G)$, the $F$ distribution is the reference distribution. This quality index can detect whether the distribution $G$ is the same as distribution $F$ in scale and/or location change (or mean change).  

In the paper, Section 2 introduces the essential background for multivariate two-sample tests and three new proposed test statistics. Section 3 make the simulations for new proposed test statistics and make comparisons with two depth-based tests under three scenarios: scale shift, location shift, and both scale and location shift. Then extend the two-sample test to a more generalized multi-sample test in Section 4 and make the comparisons. Finally, we apply the new test statistics to two real data in Section 5 and make conclusions and limitations in Section 6.

\section{Tests for homogeneity multivariate two samples }

\subsection{Parametric tests under normality assumption}

We assume the multivariate analysis of variance (MANOVA) \cite{Hair} model:
\begin{equation*}\label{m1}
X_{i,j}=\mu_{i}+e_{i,j}, i=1,2 ~\text{and}~j=1,2,\ldots,n_i,
\end{equation*}
where $e_{i,j}$ are independent $p$-dimensional normal  variables with mean $\bm{0}$ and covariance matrix $\Sigma$ denoted as $N_p(\bm{0},\Sigma)$.

We decompose the total sum of squares $S_T=S_B+S_W$, where $S_T=\sum_{i=1}^2\sum_{j=1}^{n_i}(X_{i,j}-\bar X)(X_{i,j}-\bar X)'$, $\bar X$ is the overall sample mean, $S_B=\sum_{i=1}^2\sum_{j=1}^{n_i}(\bar X_{i}-\bar X)(\bar X_{i}-\bar X)'$,   $\bar X_{i}$ is the $i$-th sample mean, and $S_W=\sum_{i=1}^2\sum_{j=1}^{n_i}(X_{i,j}-\bar X_{i})(X_{i,j}-\bar X_{i})'$. Denote the eigenvalues of $S_W^{-1}S_B$ as $\lambda_1,\ldots,\lambda_p$. 

The Wilks's test is $W=\prod_{i=1}^p1/(1+\lambda_i)$. 
The asymptotic F distribution is $(\frac{1-W}{W})( \frac{n_1+n_2-p-1}{p} )\sim F_{p,n_1+n_2-p-1,\alpha}$, where $\alpha$ is the significance level.

The Hotelling's test is $H=\sum_{i=1}^p\lambda_i$.
The asymptotic F distribution is $(\frac{n_1+n_2-p-1}{p})H \sim F_{p,n_1+n_2-p-1}$.

The Pillai's test is $T=\sum_{i=1}^n\lambda_i/(1+\lambda_i)$.
The asymptotic F distribution is $(\frac{n_1+n_2-p-1}{p})(\frac{T}{1-T}) \sim F_{p,n_1+n_2-p-1}$.

\subsection{Non-parametric tests}
In some cases, we cannot find a closed form of distribution of test statistics. Therefore, non-parametric tests are essential and there are some typical non-parametric tests.

\subsubsection{Cram\'er test}

The Cram\'er test \cite{Anderson1962, kim2020} is based on ${w_n}^2=\int_{+\infty}^{-\infty} [F_n(x)-F(x)]^2 dF(x)$, where $F_n(x)$ is empirical distribution, and $F(x)$ is theoretical distribution. Assume there are $n$ independent and identically distributed random variables $X=X_1, X_2, ..., X_n$, each have a continuous distribution function $F_x$, and $m$ independent and identically distributed random variables $Y=Y_1, Y_2, ..., Y_m$ with distribution $F_y$. $X$ and $Y$ are two mutually independent samples. In hypothesis test, we assume these two samples are from the same distribution, i.e. $H_0: F_x=F_y$. The alternative hypothesis test is $H_{\alpha}: F_x \neq F_y$. For the univariate two-sample test, the test statistics is $T=\frac{mn}{m+n} \int_{+\infty}^{-\infty} (\hat{F_x(t)}-\hat{F_y(t)})^2 d\hat{H(t)} $, where the $\hat{F_x(t)}$ and $\hat{F_y(t)}$ are empirical distribution of $X$ and $Y$ and the notation $\hat{H(t)}$ is $(n\hat{F_x(t)}+m\hat{F_y(t)})/(m+n)$. With significance level $\alpha$, the hypothesis $H_0$ is rejected if $T \ge T_{\alpha}$, and $T_{\alpha}$ is the upper $\alpha$ quantile for this distribution. i.e., $P(T < T_\alpha)=1-\alpha$.

\subsubsection{Energy Distance test}
 
The energy distance test shows the statistical distance between distributions. Assume X and X' are independent random vectors with Cumulative distribution function F, Y and Y' are independent random vectors with Cumulative distribution function G, the energy distance \cite{kim2020, szekely2013} in Euclidean space is $D^2=2 \mathbb{E}||X-Y||-\mathbb{E}||X-X'||-\mathbb{E}||Y-Y'||$. More generally, in metric spaces, the energy distance is defined as $D^2=2\mathbb{E}[d(X,Y)]-\mathbb{E}[d(X,X')]-\mathbb{E}[d(Y,Y')]$, where $d(X,Y)$ denote the distance in any metric space. To test whether the two random variables X and Y are from the same distributions, sample n samples $x_1, ..., x_n$ from X and m samples $y_1, ..., y_m$ from Y, respectively. Then the test statistics under the null hypothesis is $T=\frac{nm}{n+m} E(X,Y)$, where $E(X,Y)$ is the energy distance $D^2$. The alternative test statistics is as follows: $H=\frac{2\mathbb{E}||X-Y||-\mathbb{E}||X-X'||-\mathbb{E}||Y-Y'||}{2\mathbb{E}||X-Y||}$, which normalize the energy distance Statistics and $0\leq H \leq 1$. When $H=0$, X and Y are identically distributed. Then use $P(H<T)=1-\alpha$ to test the energy distance test with significance level $\alpha$.

\subsubsection{Depth-based Tests}

\begin{enumerate}
\item Maximum and Minimum Statistics

Let $F(x)$ be a distribution in d-dimensional space, the depth function $D(x;F)$ measures the centrality of one point $x$ in distribution $F(x)$ and transforms from space $R^d$ into [0,1]. Q Statistics, proposed by \cite{liu1993quality}, is defined as the $$Q(F,G)=P\{D(X;F)\leq D(Y;F)|X\sim F, Y\sim G\},$$ where F is the reference distribution. Suppose the two distributions $F$ and $G$ are unknown, using the empirical distributions $F_m$ and $G_n$ for $F$ and $G$, respectively. In that case, the Q statistics can be estimated as $$Q(F_m, G_n)=\frac{1}{n} \sum_{i=1}^{n} R(y_i;F_m),$$ with the sample proportion $R(y_i;F_m)$ satisfying $D(x_j,F_m) \leq D(y_i,F_m) $.
Under the null hypothesis, $H_0: F=G$, the $Q(F,G)=\frac{1}{2}$. Now we consider both $Q(F_m,G_n)$ and $Q(G_n,F_m)$ to be involved in the following test Statistics. Note that because of different reference distributions, $Q(F_m,G_n) \neq Q(G_n,F_m)$.

Based on the quality index Q, \cite{Shi2023} introduced Maximum Statistics to efficiently capture the disparity of two distributions, defined as
\begin{equation}\label{Msta}
   M_{m,n}=\left[\frac{1}{12}(\frac{1}{m}+\frac{1}{n})\right]^{-1} \max\{(Q( {F}_{m}, {G}_{n})-\frac{1}{2})^{2},
   (Q( {G}_{n}, {F}_{m})-\frac{1}{2})^{2}\}.
\end{equation}

In a similar way, we proposed the Minimum statistics, inspired by \cite{Liu04}. Our Minimum Statistics is 
\begin{equation}
M_{m,n}^*=\left[\frac{1}{12}(\frac{1}{m}+\frac{1}{n})\right]^{-\frac{1}{2}} (\frac{1}{2}-\min(Q(F_m, G_n), Q(G_n, F_m)) ) 
\end{equation}

Under null hypothesis, $F=G$, the maximum statistics $M_{m,n} \xrightarrow d \chi^2_1$, proved by \cite{Shi2023}. Similarly, under null hypothesis, our minimum statistics $M_{m,n}^* \xrightarrow d |\mathcal{N}(0,1)|$, shown in Appendix. 
The test Statistics is conducted by using $P(M_{m,n} < T_\alpha)=1-\alpha$ and $P(M_{m,n}^* < T_\alpha)=1-\alpha$, where $T_\alpha$ denotes the upper $\alpha$ quantile for Maximum and minimum distribution respectively.

\item Product and Sum Statistics

With the idea of using both quality indexes, instead of capturing the maximum or minimum of two quality indexes, we proposed Product and Sum Statistics. Product Statistics, denoted as $P_{m,n}$, defined as 
\begin{equation}
P_{m,n}=Q(F_m, G_n) Q(G_n, F_m).
\end{equation}

Similarly, Sum Statistics, denoted as $S_{m,n}$, is
\begin{equation}
S_{m,n}=Q(F_m, G_n)+Q(G_n, F_m).
\end{equation}

The test Statistics for Product $P_{m,n}$ and Sum Statistics is $P(P_{m,n} < P_\alpha)=\alpha$ and $P(S_{m,n} < S_\alpha)=\alpha$, where $P_\alpha$ is the lower $\alpha$ quantile for $P_{m,n}$ and $S_\alpha$ is the lower $\alpha$ quantile for $S_{m,n}$ respectively. These two tests are powerful as they take into account both quality indexes and captures the disparity between two samples effectively.

\item Depth-based Rank Test

Depth-based rank (DbR) test \cite{Small2011} is used to order the samples in an increasing or decreasing order to clarify the data. For the univariate sample $X_1, ..., X_n$, the rank of point $X_i$ is defined as $R(X_i)=\# \{X_j: X_j \ge X_i\}$. The notation $\#$ means the cardinality in the set. Similarly, for multivariate samples $X_1, X_2, ..., X_n$, given depth function $D(X,F)$, meaning the depth of sample X with the reference distribution F, the rank of $X_i$ is $R(X_i)=\# \{j: D(X_j,F) \ge D(X_i, F) \}$. 

Assume $X_1$ and $X_2$ are two samples with empirical distribution $\hat{F_1}$ and $\hat{F_2}$ with sample size $n_1$ and $n_2$ respectively. Then denote $R_{i,j}(k)$ as the depth rank of $X_{i,j}$ with respect to empirical distribution $\hat{F_k}$, where k=1,2. Under the null hypothesis $H_0: F_1=F_2$ and alternative hypothesis $H_\alpha: F_1 \neq F_2$, the test statistics H can be written in the form $$ H= \frac{12}{n(n+1)t} \sum_{k=1}^{2} \sum_{j=1}^{2} \frac{R_{\cdot  j}^2(k)}{n_j}-3(n+1),$$ where $t=2$, $n=n_1+n_2$, $j=1,2$ (the number of samples $X_j$), and $R_{\cdot j}(k)=\sum_{i=1}^{n_j} R_{i,j}$.

The test Statistics is conducted by using $P(H<T)=1-\alpha$, where $T$ denotes the upper $\alpha$ quantile for this distribution.

\item Modified Depth-based Rank Test

In the paper by Barale and Shirke \cite{Barale&Shirke}, they proposed the modified two-sample rank test for scale-location problems. In the univariate case, consider $(X_1, X_2, ..., X_n)$ and $(Y_1, Y_2, ..., Y_m)$ are the two samples following F and G distributions, respectively. The goal is to test whether the two samples are from the same distribution, i.e., $H_0: F=G$. Reorder the combined samples with size $N=n+m$ and rank them as $R_{(1)} < ... < R_{(n)}$ and $Q_{(1)} < ... < Q_{(m)}$ for samples from $X_i$ and $Y_j$. 
The test statistics is $B=\frac{1}{2} (B_1+B_2)$, where $B_1=\frac{1}{n} \sum_{i=1}^{n} \frac{ (R_{(i)} -\frac{N}{n} i)^2}{ \frac{i}{n+1} (1-\frac{i}{n+1}) \frac{mN}{n} }$, and $ B_2=\frac{1}{m} \sum_{j=1}^{m} \frac{ (Q_{(j)} -\frac{N}{m} j)^2}{ \frac{j}{m+1} (1-\frac{j}{m+1}) \frac{nN}{m} }$. 
The test Statistics is conducted by using $P(B<T)=1-\alpha$, where $T$ denotes the upper $\alpha$ quantile for this distribution. A larger value of test statistics will cause the rejection of the null hypothesis. However, this test is modified by Murakami to prevent the problem of not being invariant. Use B* to denote the new test statistics $B^*=\frac{1}{2} (B^*_1+B^*_2)$, 
\begin{equation*}
B^*_1=\frac{1}{n} \sum_{i=1}^{n} \frac{ (R_{(i)}-E(R_{(i)}))^2} {Var (R_{(i)})} ,
\end{equation*}
\begin{equation*}
B^*_2=\frac{1}{m} \sum_{j=1}^{m} \frac{ (Q_{(j)}-E(Q_{(j)}))^2} {Var (Q_{(j)})} , 
\end{equation*}
where 
$ E(R_{(i)})=\frac{N+1}{n+1} i $, $ E(Q_{(j)})=\frac{N+1}{m+1} j$, $Var (R_{(i)})=\frac{i}{n+1} (1-\frac{i}{n+1})\frac{m(N+1)}{n+2}$, and $Var (Q_{(j)}) =\frac{j}{m+1} (1-\frac{j}{m+1})\frac{n(N+1)}{m+2}$. Similarly, the test Statistics is conducted by using $P(B^*<T)=1-\alpha$, where $T$ denotes the upper $\alpha$ quantile for this distribution, and large test statistic values $B^*$ will cause the rejection of the null hypothesis.

In the multivariate case, let $X_{ij} \in  \mathbb{R} $, we have $X_1=\{ X_{11}, ..., X_{1n_1} \}$ and $X_2 =\{ X_{21}, ..., X_{2n_2} \}$ be the two samples from distribution $F_1$ and $F_2$ respectively. Let $\hat{F_1}$ and $\hat{F_2}$ be the empirical distributions. Assume $H_0: F_1=F_2$, here both $F_1$ and $F_2$ have equal location vector ($\mu$) and equal scale matrix ($\Sigma$), and $H_\alpha: F_1 \neq F_2$.
The proposed procedure for this test statistics is as follows: 

First, combine the two samples $X_1$ and $X_2$ with total size N, i.e. $N=n_1+n_2$, denote the combined samples as $Z=X_1 \bigcup X_2$, and $Z_t$ is each observation in $Z$ with $t=1, \dots, N$.
Then, compute the depth of all $Z_t$ with respect to $\hat{F_1}$ and $\hat{F_2}$, and denote them as $D(Z_t \hat{F_1})$ and $D(Z_t, \hat{F_2})$ respectively. 
Then, rank all the observations based on depth values $D(Z_t, \hat{F_1})$ and record as $R^{\hat{F_1}}_t$, same to depth values $D(Z_t, \hat{F_2})$ and record as $R^{\hat{F_2}}_t$.
Finally, reorder these ranks, select those ranks corresponding to sample $X_2$ and record as $R^{\hat{F_1}}_{(1)} < ...< R^{\hat{F_1}}_{(n_2)}$, select those ranks corresponding to sample $X_1$ and record as $R^{\hat{F_2}}_{(1)} < ...< R^{\hat{F_2}}_{(n_1)}$.

Then the test statistics $B=max(B^{\hat{F_1}}, B^{\hat{F_2}})$. The test Statistics is conducted by using $P(B<T)=1-\alpha$, where $T$ denotes the upper $\alpha$ quantile for this distribution.

In detail, 
\begin{equation*}
B^{\hat{F_1}}=\frac{1}{n_2} \sum_{j=1}^{n_2} \frac{(R^{\hat{F_1}}_{(j)} -E(R^{\hat{F_1}}_{(j)} ))^2}{Var (R^{\hat{F_1}}_{(j)} )},
\end{equation*}
\begin{equation*}
B^{\hat{F_2}}=\frac{1}{n_1} \sum_{j=1}^{n_1} \frac{(R^{\hat{F_2}}_{(j)} -E(R^{\hat{F_2}}_{(j)} ))^2}{Var (R^{\hat{F_2}}_{(j)} )}. 
\end{equation*}

Here, $E(R^{\hat{F_1}}_{(j)})=\frac{N+1}{n_2+1} j$, $Var(R^{\hat{F_1}}_{(j)})=\frac{j}{n_2+1}(1-\frac{j}{n_2+1})\frac{n_1(N+1)}{n_2+2}$, $E({R^{\hat{F_2}}_{(j)}})=\frac{N+1}{n_1+1} j$, and $Var(R^{\hat{F_2}}_{(j)})=\frac{j}{n_1+1}(1-\frac{j}{n_1+1})\frac{n_2(N+1)}{n_1+2}$.

\end{enumerate}

\section{Simulation studies: two-sample test}
We proposed the distribution of Minimum Statistic, Product Statistic and Sum Statistic for the two-sample cases. Simulations are conducted to see the performance of these statistics compared with other depth-based tests. Assume we have random samples $x_1, x_2, \dots, x_m$ and $y_1, y_2, \dots, y_n$ from distributions $F$ and $G$ respectively with size m, and n. 

\begin{figure} 
\begin{center}
\includegraphics[width=\textwidth]{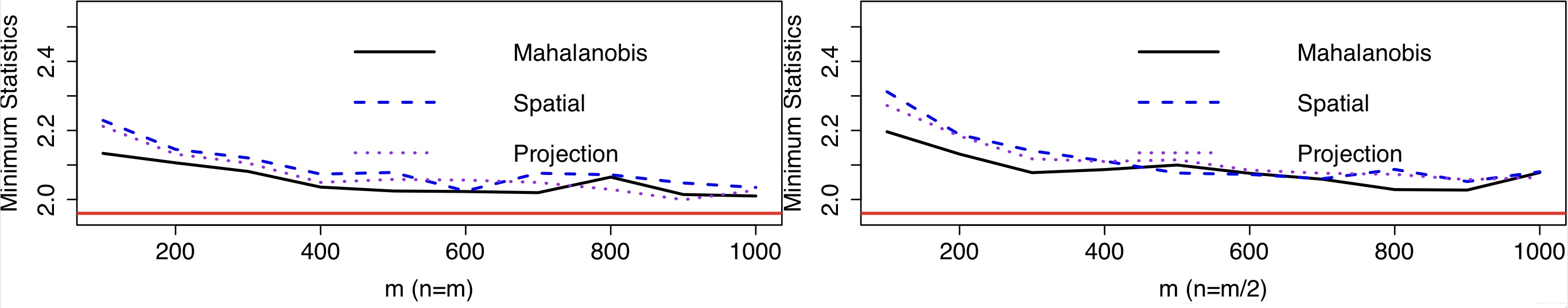}
\caption{Comparison of empirical upper 95\% quantiles of minimum statistics $M_n$ for $m=100,200,\ldots,1000$ and $n=m$ (1st column) or $n=m/2$ (2nd column).}
\label{fig:type1}
\end{center}
\end{figure}

First, assume we have two equal distributions $F=G=N(\bm{0},I_{2\times2})$, where $N(\bm{0},I_{2\times2})$ represents the bivariate normal distribution with mean vector $\bm{0}$ and two-by-two identity covariance matrix. By setting up the sample size $m=100,200,\ldots,1000$ and $n$ with $n=m$ or $n=m/2$. Since we proved that the Minimum statistics $M_n$ follows a half-normal asymptotic null distribution, the upper 95\% quantile is 1.96 in this case. We plotted the Type I error curve of Minimum Statistic in Figure  \ref{fig:type1}, presenting the empirical quantiles based on different values of $m, n$ and different depth functions, comparing with theoretical quantiles. Here we simulated for 10000 repetitions. This figure shows the convergence rate of three depth functions and both turned out that the Mahalanobis depth converges fastest relative to Spatial and Projection depth.

To show the effect of our proposed Minimum Statistic $M_{m,n}^*$, Product Statistic $P_{m,n}$, and Sum Statistic $S_{m,n}$, we made comparisons with Maximum Statistic $M_{m,n}$, depth-based rank (DbR) statistic \citep{Small2011}, and the Modified Depth-based Rank Statistics (BDbR)\cite{Barale&Shirke}. Assume we use $\alpha=0.05$ to do the simulation on power of these tests. The critical values of $M_{m,n}$, $M_{m,n}^*$, DbR, and BDbR are based on upper 95\% quantiles, and $P_{m,n}$ and $S_{m,n}$ are based on lower 5\% quantile. Similar to Type I error, we compare the power for three depth functions Mahalanobis depth, spatial depth, and projection depth at different sample sizes ($m=n=100,200,\ldots,1000$, or $m=100,200,\ldots,1000, n=m/2$).

We compare the power based on three scenarios: scale change, mean change, and both scale and mean change.

\begin{figure} 
\begin{center}
\includegraphics[width=\textwidth]{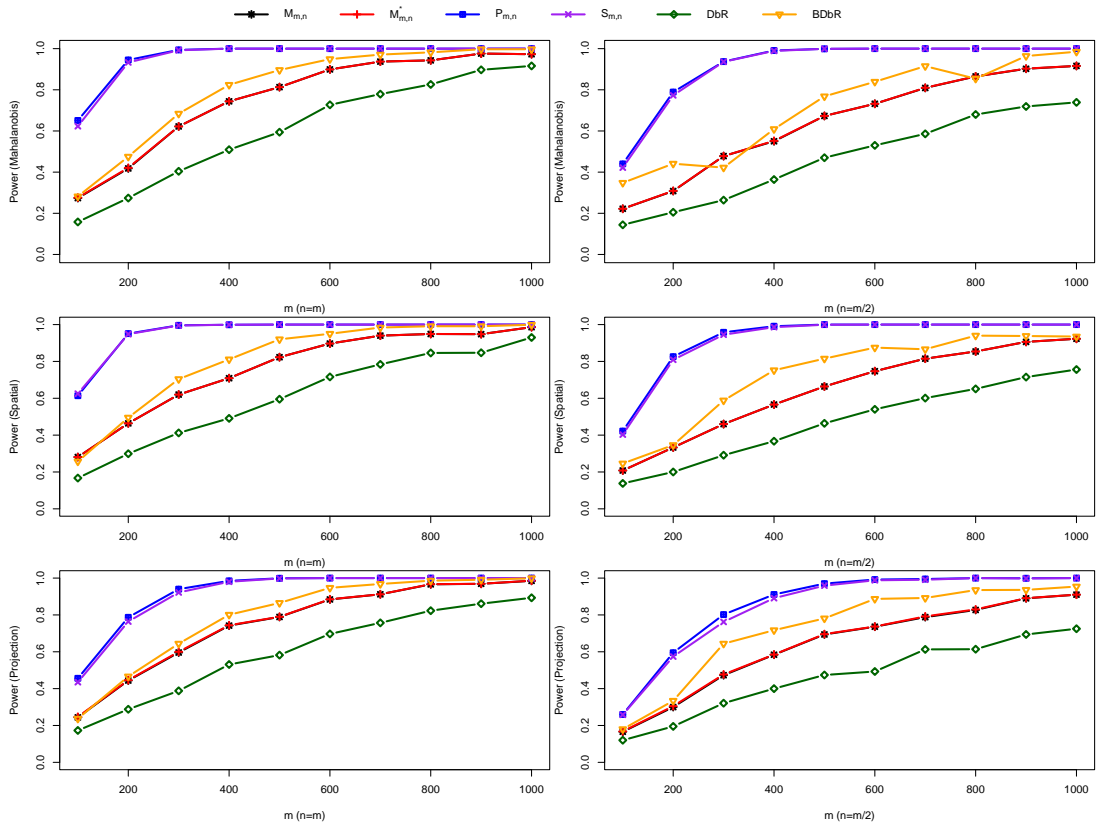}
\caption{Power comparison under alternative hypothesis $F=N(\bm{0},I_{2\times2})$ against $G=N(\bm{0},I_{2\times2}+0.5\tilde{I}_{2\times2})$ for $m=100,200,\ldots,1000$ and $n=m$ (1st column) or $n=m/2$ (2nd column) for Mahalanobis depth (Row 1), Spatial depth (Row 2), and Projection depth (Row 3). } 
\label{fig:power1}
\end{center}
\end{figure}

(1)  Two bivariate normal distributions with a scale change:

Assume one sample is from $F=N(\bm{0},I_{2\times2})$ and another sample from $G=N(\bm{0},I_{2\times2}+0.5\tilde{I}_{2\times2})$, where $\tilde{I}_{2\times2}=((0,1)^\top,(1,0)^\top)$. The power comparisons are shown in Figure \ref{fig:power1} for different depth functions and sample sizes for 1000 repetitions. Each row shows the power based on Mahalanobis depth, spatial depth, and projection depth, respectively. We can clearly observe that the three depth functions perform similar trends for powers of different test statistics. Under all three depth functions, the Product Statistic and Sum Statistic outperform all other statistics, and the powers achieve one quickly compared with other statistics. The product and sum work best in power because it takes into account both quality indexes and captures the disparity between two samples effectively. In addition, the Minimum statistics is comparable to the BDbR.

\begin{figure}
\begin{center}
\includegraphics[width=\textwidth]{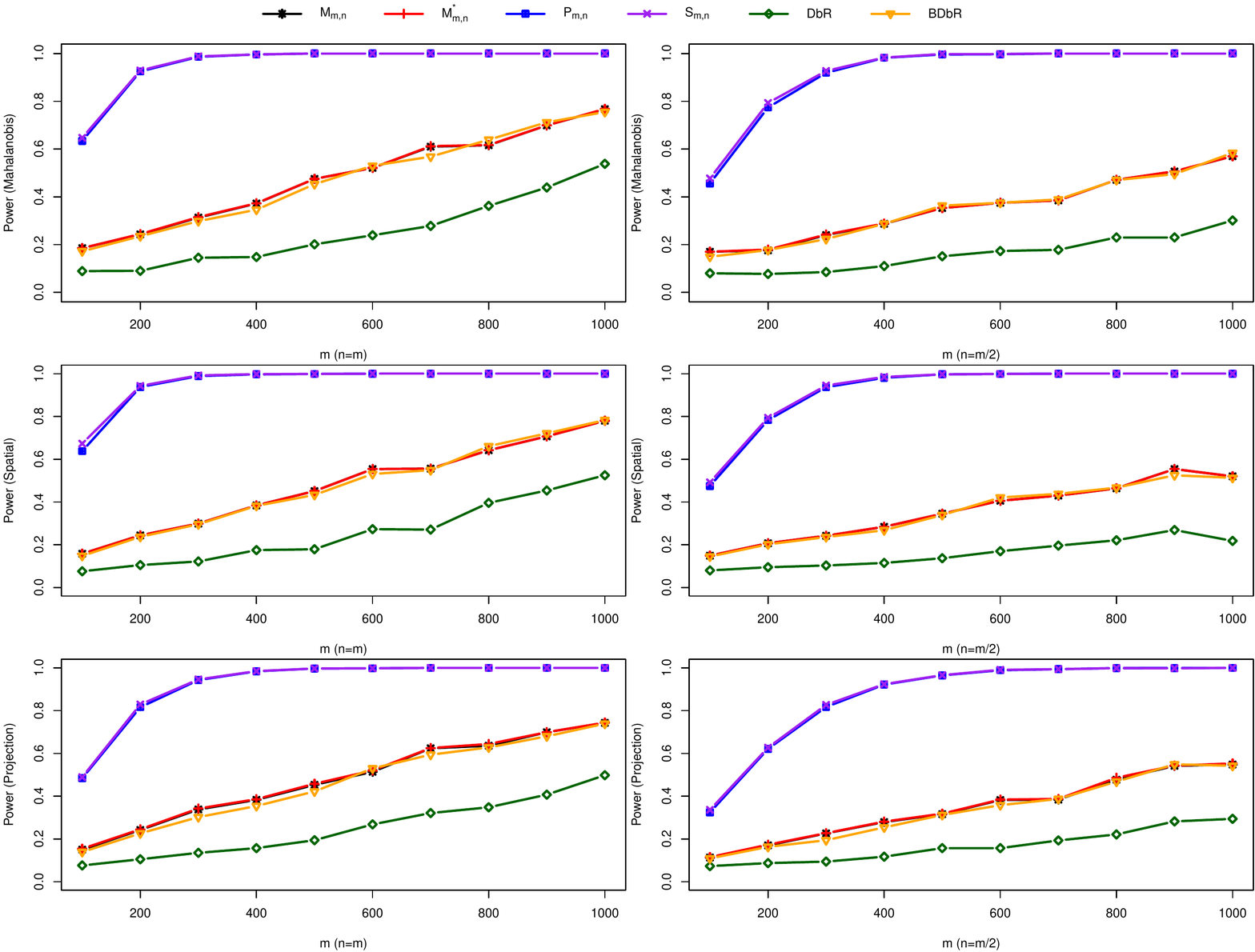}
\caption{Power comparison under alternative hypothesis $F=N( \bm{0},I_{2\times2})$ against $G=N((0.3,0.3)^\top, I_{2\times2})$ for $m=100,200,\ldots,1000$ and $n=m$ (1st column) or $n=m/2$ (2nd column) for Mahalanobis depth (Row 1), Spatial depth (Row 2), and Projection depth (Row 3).}
\label{fig:power2}
\end{center}
\end{figure}

(2) Two bivariate normal distributions with a mean change:

To visualize how the mean change will affect the power, we assume one sample from $F=N(\bm{0},I_{2\times2})$ and another from $G=N((0.3,0.3)^\top, I_{2\times2})$. Similarly, as shown in Figure \ref{fig:power2}, the Product Statistic and Sum Statistic have the largest power and are almost the same at any sample size. In this case, all other statistics have relatively low power, and Maximum Statistic, Minimum statistic, and BDbR are almost on the same line.

\begin{figure}
\begin{center}
\includegraphics[width=\textwidth]{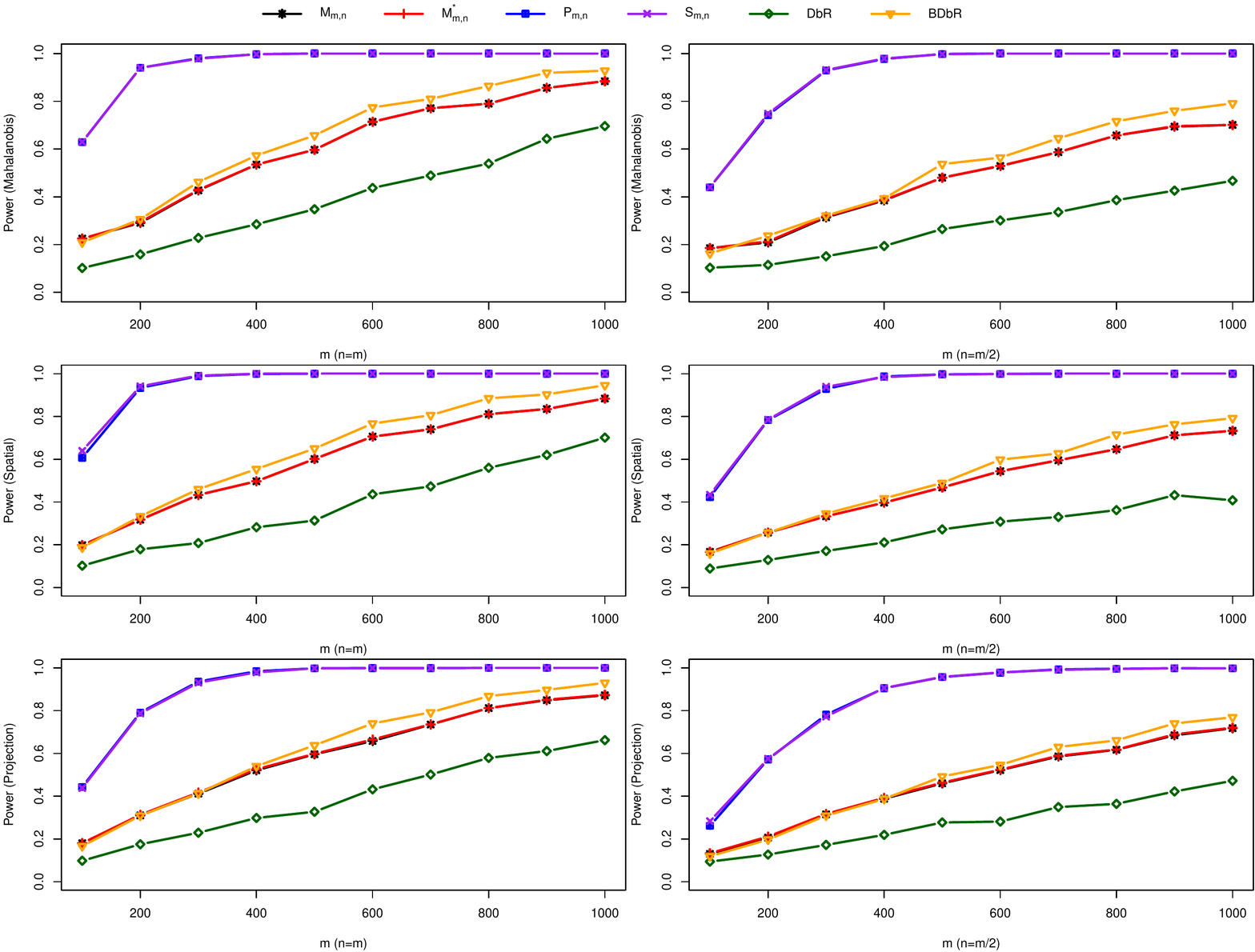}
\caption{Power comparison under alternative hypothesis $F=N( \bm{0},I_{2\times2})$ against $G=N((0.2,0.2)^\top, I_{2\times2}+0.4\tilde{I}_{2\times2})$ for $m=100,200,\ldots,1000$ and $n=m$ (1st column) or $n=m/2$ (2nd column) for Mahalanobis depth (Row 1), Spatial depth (Row 2), and Projection depth (Row 3).}
\label{fig:power3}
\end{center}
\end{figure}

(3) Two bivariate normal distributions with both scale and mean change:

Under both mean and scale change, as one sample from $F=N(\bm{0},I_{2\times2})$ and another sample from $G=N((0.2,0.2)^\top, I_{2\times2}+0.4\tilde{I}_{2\times2})$, the result in Figure \ref{fig:power3} is similar as the previous scenarios. The Product Statistic and Sum Statistic outperform all other statistics and the Minimum statistics is comparable to the BDbR.

\section{Multi-sample Test}
Previous section is simulations on two-sample tests, now extend the two-sample cases to Multi-sample cases. 

We generalize the Minimum Statistics for k-sample in the form: 

\begin{equation}
M_{n_1, ... ,n_k}^*=\max_{1\leq i,j\leq k, i\neq j} \left[\frac{1}{12} (\frac{1}{n_i}+\frac{1}{n_j})\right]^{-\frac{1}{2}} (\frac{1}{2} -Q(F^{(i)}_{n_i},F^{(j)}_{n_j}))
\end{equation}

This can be written in this form, based on the proof of asymptotic distribution of Minimum Statistics: 
\begin{equation}
P(M_{n_1,\ldots,n_k}^* \leq x)\rightarrow P\left\{\max_{1\leq i<j\leq k}(c_{i,j}Z_i+\tilde{c}_{i,j}Z_j)\leq x\right\},
\end{equation}
where $Z_1,Z_2,\ldots,Z_k$  are independent from  $N(0,1)$, $c_{i,j}=\lim n_i^{-1/2}(n_i^{-1}+n_j^{-1})^{-1/2}$, and $\tilde{c}_{i,j}=\lim n_j^{-1/2}(n_i^{-1}+n_j^{-1})^{-1/2}$ with $c_{i,j}^2+\tilde{c}_{i,j}^2=1$.

In this section, we will do the power comparisons for three-sample cases. Hence, when $k=3$, the Minimum Statistics can be expanded as 

\begin{equation*}
\begin{split}
M_{n_1,n_2,n_3}^*=\max \{ [\frac{1}{12}(\frac{1}{n_{1}}+\frac{1}{n_{2}})]^{-\frac{1}{2}}[\frac{1}{2}-Q( {F}^{(1)}_{n_{1}},{F}^{(2)}_{n_{2}})], [(\frac{1}{12}(\frac{1}{n_{1}}+\frac{1}{n_{3}})]^{-\frac{1}{2}}[\frac{1}{2}-Q( {F}^{(1)}_{n_{1}}, {F}^{(3)}_{n_{3}})], \\
[(\frac{1}{12}(\frac{1}{n_{2}}+\frac{1}{n_{1}})]^{-\frac{1}{2}}[\frac{1}{2}-Q( {F}^{(2)}_{n_{2}}, {F}^{(1)}_{n_{1}})],
[(\frac{1}{12}(\frac{1}{n_{2}}+\frac{1}{n_{3}})]^{-\frac{1}{2}}[\frac{1}{2}-Q( {F}^{(2)}_{n_{2}}, {F}^{(3)}_{n_{3}})], \\
[(\frac{1}{12}(\frac{1}{n_{3}}+\frac{1}{n_{2}})]^{-\frac{1}{2}}[\frac{1}{2}-Q( {F}^{(3)}_{n_{3}}, {F}^{(2)}_{n_{2}})],
[(\frac{1}{12}(\frac{1}{n_{3}}+\frac{1}{n_{1}})]^{-\frac{1}{2}}[\frac{1}{2}-Q( {F}^{(3)}_{n_{3}}, {F}^{(1)}_{n_{1}})]
\}.
\end{split}
\end{equation*}
and

\begin{align}\label{p1}
P(M_{n_1,n_2,n_3}^* \leq x) \rightarrow P\{ &-x \leq (c_{1,2}Z_1+\tilde{c}_{1,2}Z_2)\leq x, -x \leq (c_{1,3}Z_1+\tilde{c}_{1,3}Z_3)\leq x, \\
&-x\leq (c_{2,3}Z_2+\tilde{c}_{2,3}Z_3)\leq x  \}.
\end{align}

Similarly, we computed the formula of Product Statistics and Sum Statistics for k-sample comparison.
\begin{equation}
P_{n_1, ... ,n_k}=\prod_{1\leq i,j\leq k, i\neq j}^{k}  Q(F^{(i)}_{n_i},F^{(j)}_{n_j})
\end{equation}
and

\begin{equation}
S_{n_1, ... ,n_k}=\sum_{1\leq i,j\leq k, i\neq j}^{k}  Q(F^{(i)}_{n_i},F^{(j)}_{n_j})
\end{equation}

For three-sample cases $k=3$,  the product and sum statistics are
\begin{equation*}
\begin{split}
P_{n_1,n_2,n_3}= Q( {F}^{(1)}_{n_{1}},{F}^{(2)}_{n_{2}})   Q( {F}^{(1)}_{n_{1}}, {F}^{(3)}_{n_{3}}) 
 Q( {F}^{(2)}_{n_{2}}, {F}^{(1)}_{n_{1}})   \\
 Q( {F}^{(2)}_{n_{2}}, {F}^{(3)}_{n_{3}}) 
 Q( {F}^{(3)}_{n_{3}}, {F}^{(2)}_{n_{2}})   Q( {F}^{(3)}_{n_{3}}, {F}^{(1)}_{n_{1}}) 
\end{split}
\end{equation*}
and

\begin{equation*}
\begin{split}
S_{n_1,n_2,n_3}= Q( {F}^{(1)}_{n_{1}},{F}^{(2)}_{n_{2}})  +Q( {F}^{(1)}_{n_{1}}, {F}^{(3)}_{n_{3}}) +
 Q( {F}^{(2)}_{n_{2}}, {F}^{(1)}_{n_{1}})  \\
 +Q( {F}^{(2)}_{n_{2}}, {F}^{(3)}_{n_{3}}) 
 +Q( {F}^{(3)}_{n_{3}}, {F}^{(2)}_{n_{2}}) +  Q( {F}^{(3)}_{n_{3}}, {F}^{(1)}_{n_{1}}) 
\end{split}
\end{equation*}

Based on the expanded three-sample cases for minimum statistics $M_{m,n}^*$, product statistics $P_{m,n}$, and sum statistics $S_{m,n}$, we can make power comparisons for these three Statistics and compare with maximum statistics $M_{m,n}$ and DbR statistic.
Similar to two-sample cases, the critical values are based on upper 95\% quantiles for $M_{m,n}$, $M_{m,n}^*$, and DbR, and lower 5\% quantile for $P_{m,n}$ and $S_{m,n}$. The power comparisons are performed with different sample sizes ($m=100,200,\ldots,1000$, $n=k=m$ or $n=2k=m/2$) with Mahalanobis depth, spatial depth, or projection depth. 

We consider the three distributions as $F_1$, $F_2$, and $F_3$, each with sample size $m, n, k$. To check the power of these statistics, we assume two cases: (1) Two same distributions and one different distribution; or (2) All three different distributions.

\begin{figure}
\begin{center}
\includegraphics[width=\textwidth]{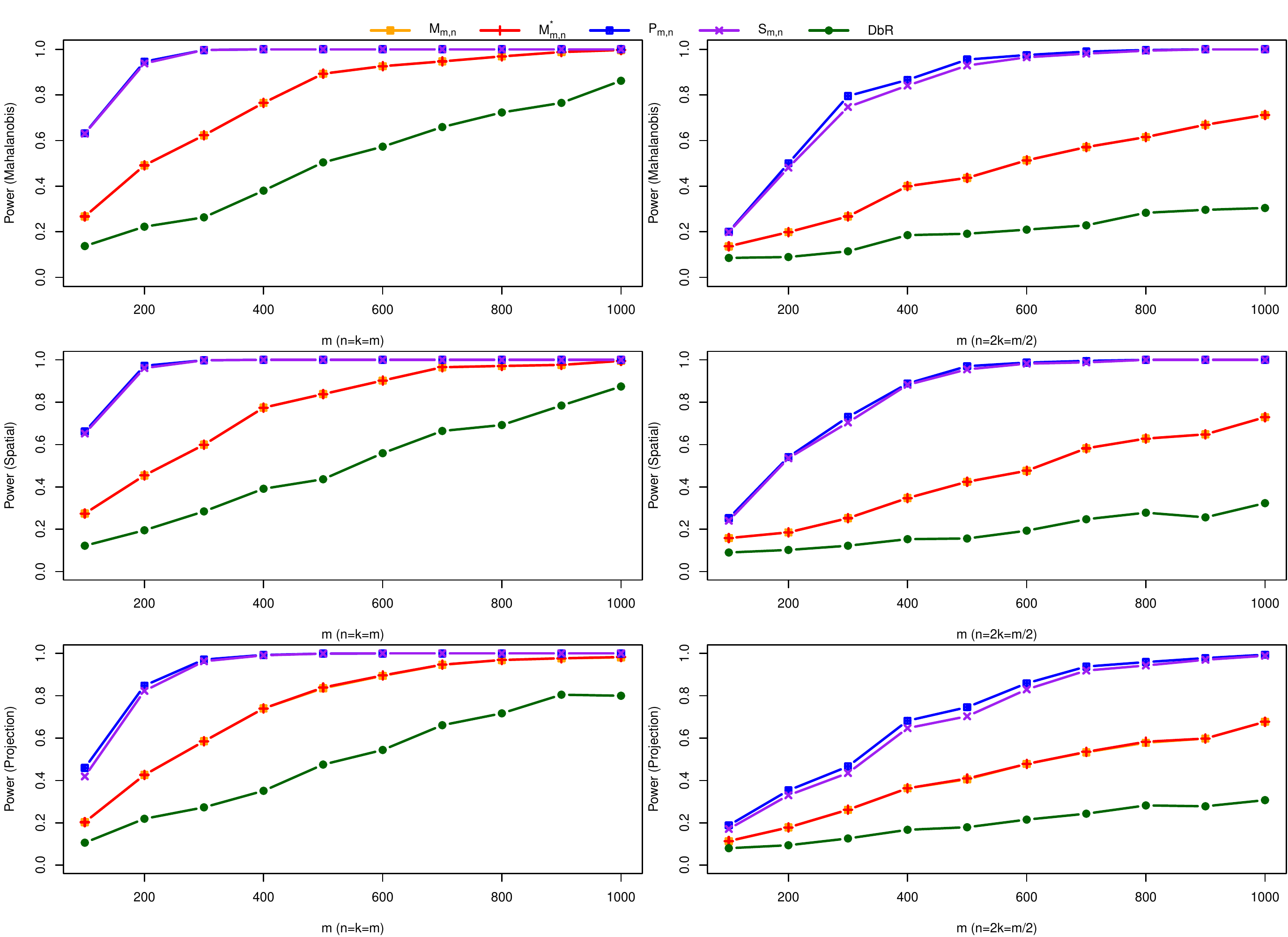}
\caption{Comparisons of power of five statistics $M_{m,n}$, $M_{m,n}^*$, $P_{m,n}$,  $S_{m,n}$and  DbR statistic in \citep{Small2011} under alternative hypothesis $F_1=F_2=N(\bm{0},I_{2\times2})$ and $F_3=N((0,0)^\top,I_{2\times2}+0.5\tilde{I}_{2\times2})$ for Mahalanobis depth (Row 1), Spatial depth (Row 2), and Projection depth (Row 3).}
\label{fig:power4}
\end{center}
\end{figure}

(1) Three bivariate normal distributions:  

Assume $F_1=F_2=N(\bm{0},I_{2\times2})$ and $F_3=N((0,0)^\top,I_{2\times2}+0.5\tilde{I}_{2\times2})$. With sample sizes from $m=100,200,\ldots,1000$, and $n=k=m$ or $n=2k=m/2$, the power of five statistics is shown in Figure \ref{fig:power4} for different depth functions. The data is simulated 1000 times. The trend is similar under all depth functions. All graphs showed that the product and sum Statistics perform the best among all five distributions with almost the same value at any sample size. Similar to minimum Statistics, with almost the same value as Maximum Statistic, and performs better than DbR. At sample size $n=k=m$, all five statistics have a larger increasing rate than at $n=2k=m/2$.

\begin{figure}
\begin{center}
\includegraphics[width=\textwidth]{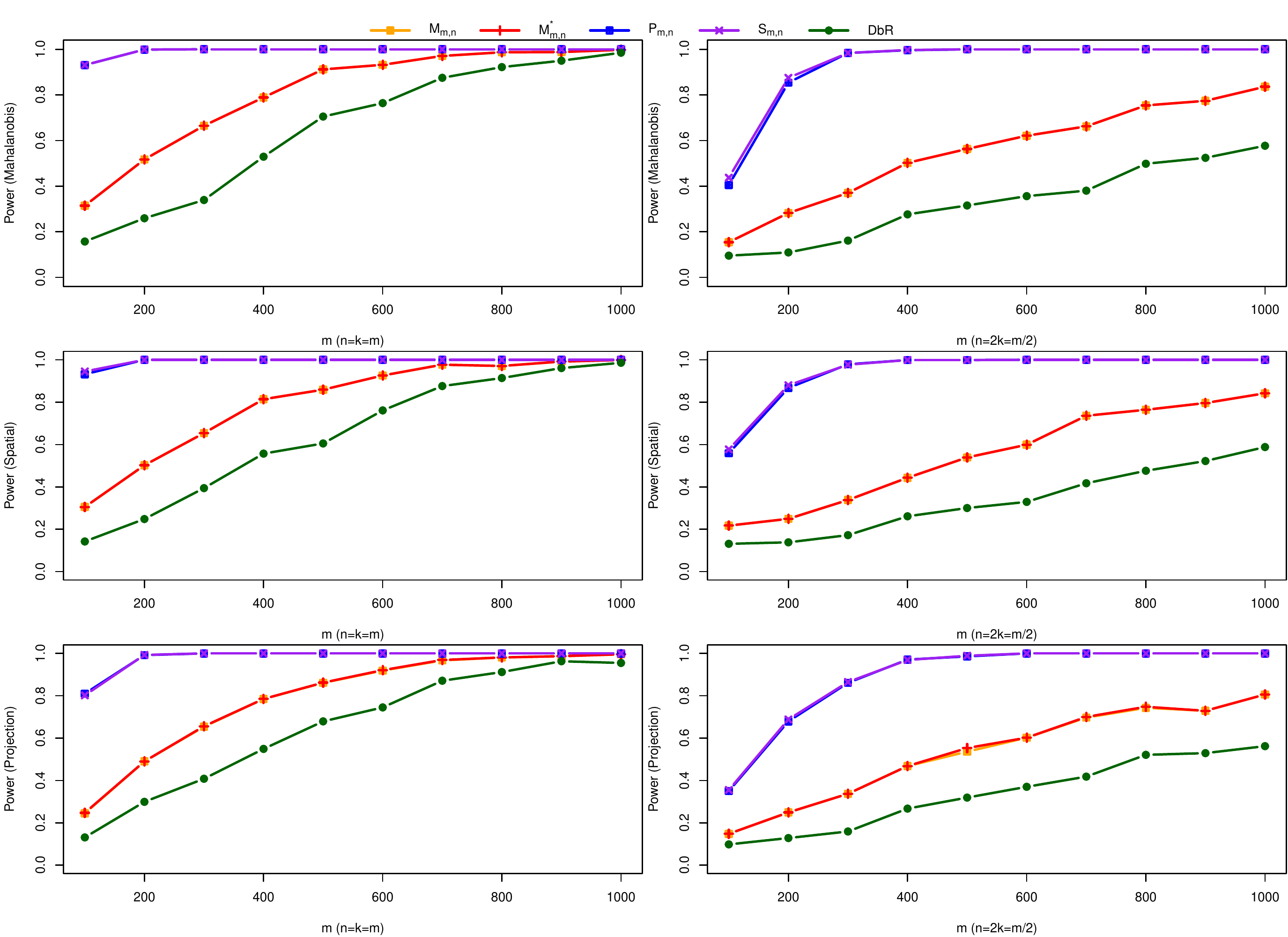}
\caption{Comparisons of power of five statistics $M_{m,n}$, $M_{m,n}^*$, $P_{m,n}$,  $S_{m,n}$and  DbR statistic in \citep{Small2011} under alternative hypothesis $F_1=N(\bm{0},I_{2\times2})$, $F_2=N((0.3,0.3)^\top,I_{2\times2})$, and
$F_3=N((0,0)^\top,I_{2\times2}+0.5\tilde{I}_{2\times2})$ for Mahalanobis depth (Row 1), Spatial depth (Row 2), and Projection depth (Row 3).}
\label{fig:power5}
\end{center}
\end{figure}

(2)  Three distinguished bivariate normal distributions:  

Let $F_1=N(\bm{0},I_{2\times2})$, $F_2=N((0.3,0.3)^\top,I_{2\times2})$, and $F_3=N((0,0)^\top,I_{2\times2}+0.5\tilde{I}_{2\times2})$. 
Similar results are shown in Figure \ref{fig:power5}. The product and sum Statistics perform the best among all five distributions at any sample size. At sample size $n=k=m$, these two statistics have power almost close to one at small sample sizes. At sample size $n=k=m$, all five statistics have a larger increasing rate than at $n=2k=m/2$.

Both cases showed that our Product and Sum Statistics are better than other statistics for all three depth functions and different sample sizes. With equal sample sizes for three distributions, the power has a larger increasing rate than under different sample sizes.

\section{Real Data Analysis}
With the above simulations, we extend our proposed minimum statistic, product statistic, and sum statistic in real data cases. We used two data sets to analyze the performance of these test statistics in three-sample comparisons.

\subsection{Sloan Digital Sky Survey data}
The Sloan Digital Sky Survey data dataset is a dataset in \textit{astrodatR} in \textit{R} that contains three classes of point source, with measurements on four color indices (u-g, g-r, r-i, i-z). The three classes are classified as quasars (Class 1), main sequence and giant stars (Class 2), and giant stars (Class 3), with sample sizes 2000, 5000, and 2000 respectively. We proposed the three-sample test on this data set to see any correlation between the distribution of four color indices among three classes.

Scale curves, introduced by \cite{Liu99}, is a measure of dispersion to compare the scale of multiple distributions. 
$D_\alpha(F)$, the $\alpha$-trimmed region with respect to distribution $F$, is defined as the $$D_\alpha (F)=\left\{x\in\mathds{R}^d: D(x; F) \geq \alpha \right \}.$$ Then the volume of this convex region is $V(\alpha; F_m)$ and the scale curve is the volume at $1-\alpha$ scale. 

\begin{figure} 
\begin{center}
 \includegraphics[width=0.7\textwidth]{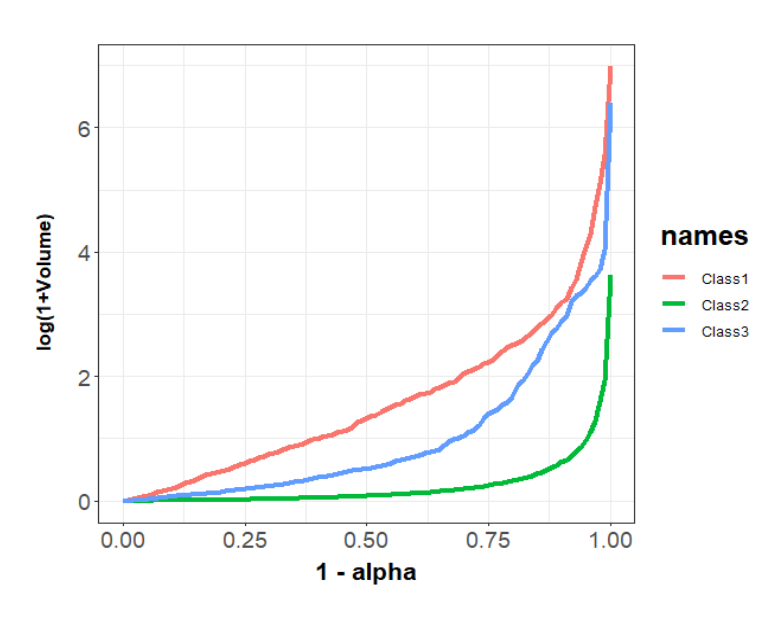}
\caption{Scale curves for three classes of under Mahalanobis depth in log scale in Sloan Digital Sky Survey data} 
 \label{fig:star}
\end{center}
\end{figure}

The scale curve of three classes in Sloan Digital Sky Survey data under Mahalanobis depth is in Figure \ref{fig:star}. To visualize the dispersion in a detailed way, we plotted the scale curve in log scale. The non-overlapping curves representing the three classes may be different. Hence, calculating the $p$-value and asymptotic $p$-value \eqref{p1} are essential.

We simulated the whole data set with 1000 repetitions for each depth function and sample sizes $n_1=2000, n_2=5000, n_3=2000$ to find the $p$-value. The results showed that all minimum statistic, product statistic, sum statistic, and DbR are zero for all Mahalanobis depth, spatial depth, and projection depth. We also calculated the asymptotic $p$-value for minimum statistics, and the result is also all zero for all depth functions. These values showed that there is a strong correlation between the combined four color indices and three classes of point source. Then if we consider the two-sample cases with only two classes: quasars (Class 1) and giant stars (Class 3), the sample sizes are all 2000. Similarly, the $p$-values and asymptotic $p$-values are also zero under all depth functions, indicating the significant differences in color indices for quasars and giant stars.

\subsection{Skull data}
In \textit{R} package, \textit{HSAUR} contains the Egyptian skulls data with four measurements in five epochs (4000 B.C., 3300 B.C., 1850 B.C., 200 B.C., and 150 A.D.), each with 30 samples. The four measurements contain maximum breaths, basibregmatic heights, basialiveolar length, and nasal heights of the skull. We are interested in whether the skulls change as time changes because of the effect of interbreeding with immigration. 

\begin{figure} 
\begin{center}
 \includegraphics[width=0.7\textwidth]{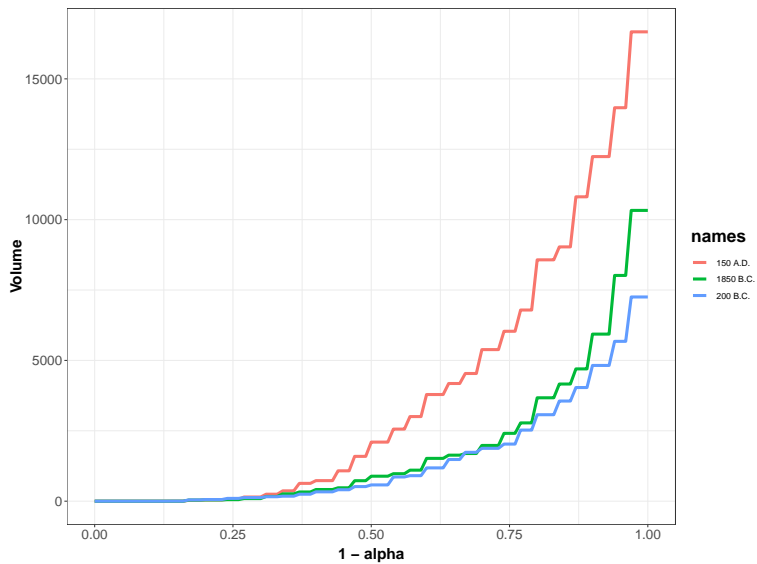}
\caption{Scale curves of skull data for epochs: 1850 B.C., 200 B.C., and 150 A.D. under Mahalanobis depth} 
 \label{fig:skull1}
\end{center}
\end{figure}

We first consider these three epochs: 1850 B.C., 200 B.C., and 150 A.D., to see if the skull size varies as time varies. The scale curve is shown in Figure \ref{fig:skull1}. The small difference between epochs 1850 B.C. and 200 B.C. represents that there may not be significant changes in skull sizes in these two time periods. Therefore, we summarized the estimated $p$-values for minimum statistic, product statistic, sum statistic, and DbR in these three epochs in Table \ref{table1}, for 5000 iterations. The estimated $p$-values are all larger than the significance level of 0.05. Hence, there is no strong correlation between skull sizes and interbreeding with immigrations for these three epochs. The asymptotic $p$-values for Minimum Statistic are 0.026864, 0.0151585, and 0.0313754 under Mahalanobis depth, Spatial depth, and Projection depth, respectively. 
 
\begin{table}
\centering
\begin{tabular}{ |p{4cm}|p{2cm}|p{2cm}|p{2cm}|p{2cm}|} 
\hline 
  & $M_{m,n}^*$ & $P_{m,n}$ & $S_{m,n}$ &DbR \\
\hline 
 Mahalanobis distance & 0.2612 & 0.2410 & 0.2886 & 0.1422\\ 
 Spatial distance & 0.2740 & 0.2110 & 0.2484 & 0.1428 \\ 
 Projection distance & 0.2446 & 0.0912 & 0.0980 & 0.1136  \\ 
\hline 
\end{tabular}
\caption{ $p$-values for  Minimum Statistic, Product Statistic, Sum Statistic, and DbR in skull data for three epochs (1850 B.C., 200 B.C., and 150 A.D.) under Mahalanobis depth, Spatial depth, and Projection depth}
\label{table1}
\end{table}

\begin{figure} 
\begin{center}
 \includegraphics[width=0.7\textwidth]{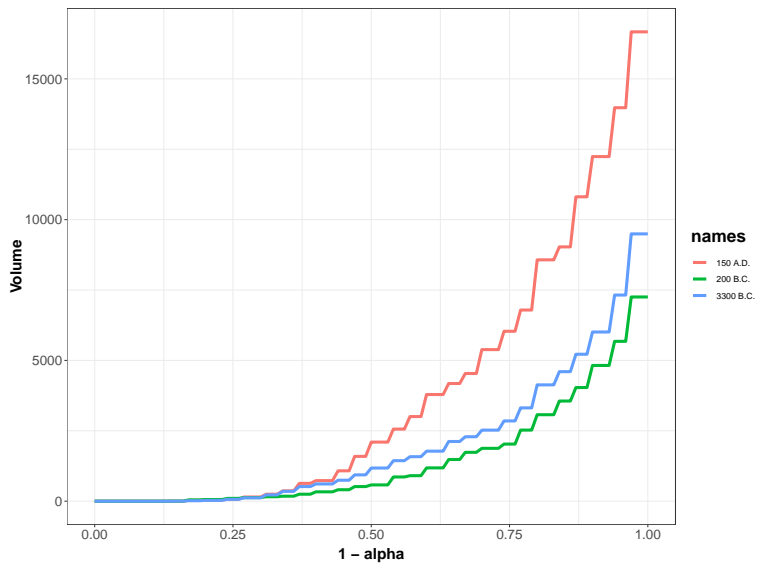}
\caption{Scale curves of skull data for epochs: 3300 B.C., 200 B.C., and 150 A.D. under Mahalanobis depth} 
 \label{fig:skull2}
\end{center}
\end{figure}

In a similar way, we did another data analysis on epochs: 3300 B.C., 200 B.C., and 150 A.D. Visualized from Figure \ref{fig:skull2}, the three curves are relatively far apart, meaning there may be a significant difference in skulls as time changes. The estimated $p$-values for minimum statistic, product statistic, sum statistic, and DbR are summarized in Table \ref{table2} for 5000 iteration times. Observed that these $p$-values are all smaller than 0.05 significance level and close to zero, concluded that there is large difference in skulls between these three epochs. We also calculated the asymptotic $p$-values for Minimum Statistic with the same conclusions, with values 0.001243, 0.0005263, and 0.0042542 under Mahalanobis depth, Spatial depth, and Projection depth, respectively.

\begin{table}
\centering
\begin{tabular}{ |p{4cm}|p{2cm}|p{2cm}|p{2cm}|p{2cm}|} 
\hline 
  & $M_{m,n}^*$ & $P_{m,n}$ & $S_{m,n}$ &DbR \\
\hline 
 Mahalanobis distance & 0.0378 & 0.0016 & 0.0018 & 0.0120\\ 
 Spatial distance & 0.0328 & 0.0012 & 0.0014 & 0.0094 \\ 
 Projection distance & 0.0260 & 0.0018 & 0.0038 & 0.0100  \\ 
\hline 
\end{tabular}
\caption{ $p$-values for  Minimum Statistic, Product Statistic, Sum Statistic, and DbR in skull data for three epochs (3300 B.C., 200 B.C., and 150 A.D.) under Mahalanobis depth, Spatial depth, and Projection depth}
\label{table2}
\end{table}

\section{Conclusions and Limitations}

This paper involves three new test statistics to test the homogeneity of multivariate two-samples based on data depth. We proved that the minimum statistics is the asymptotic half-normal distribution. In the simulation study, we compared six test statistics. We concluded that the product and sum statistics outperform all other statistics, and the powers achieve one faster than others. The minimum statistics is comparable with the BDbR proposed by Barale and Shirke. This research can be further improved by finding the asymptotic distributions of product and sum statistics with their higher-order approximations. It is challenging to find the asymptotic distributions of product and sum statistics as it involves multiple integrals for each component under Mahalanobis depth.

\section*{Acknowledgement(s)}

An unnumbered section, e.g.\ \verb"\section*{Acknowledgements}", may be used for thanks, etc.\ if required and included \emph{in the non-anonymous version} before any Notes or References.

\section*{Disclosure statement}

An unnumbered section, e.g.\ \verb"\section*{Disclosure statement}", may be used to declare any potential conflict of interest and included \emph{in the non-anonymous version} before any Notes or References, after any Acknowledgements and before any Funding information.

\section*{Funding}

An unnumbered section, e.g.\ \verb"\section*{Funding}", may be used for grant details, etc.\ if required and included \emph{in the non-anonymous version} before any Notes or References.

\section*{Notes on contributor(s)}

An unnumbered section, e.g.\ \verb"\section*{Notes on contributors}", may be included \emph{in the non-anonymous version} if required. A photograph may be added if requested.

\section*{Nomenclature/Notation}

An unnumbered section, e.g.\ \verb"\section*{Nomenclature}" (or \verb"\section*{Notation}"), may be included if required, before any Notes or References.

\section*{Notes}

An unnumbered `Notes' section may be included before the References (if using the \verb"endnotes" package, use the command \verb"\theendnotes" where the notes are to appear, instead of creating a \verb"\section*").


\appendix

\section{Proof of Minimum Statistic}

By \cite{liu1993quality,zuo2006limiting,Shi2023}, we have the property of Q statistics: $$Q( {G}_n, {F}_m)-1/2 = 1/2-Q( {F}_m, {G}_n)+o_p(n^{-1/2})+o_p(m^{-1/2}),$$ and \cite{zuo2006limiting} showed that $$ \left[\frac{1}{12}(\frac{1}{m}+\frac{1}{n})\right]^{-\frac{1}{2}} (Q(F_m, G_n)-\frac{1}{2} ) \xrightarrow d \mathcal{N}(0,1).$$

We have
\begin{align*}
&\min(Q(F_m, G_n), Q(G_n, F_m)) \\
= &\min (Q(F_m, G_n)-\frac{1}{2}, Q(G_n, F_m)-\frac{1}{2} )+\frac{1}{2}\\
=& -\left| Q(F_m, G_n)-\frac{1}{2} \right | +\frac{1}{2} 
\end{align*}

Hence,

\begin{align*}
M_n  & = \left[\frac{1}{12}(\frac{1}{m}+\frac{1}{n})\right]^{-\frac{1}{2}} (\frac{1}{2}- \left| Q(F_m, G_n)-\frac{1}{2} \right | - \frac{1}{2}    ) \\
&=\left[\frac{1}{12}(\frac{1}{m}+\frac{1}{n})\right]^{-\frac{1}{2}}  \left| Q(F_m, G_n)-\frac{1}{2} \right |  \xrightarrow d \left | \mathcal{N}(0,1) \right |
\end{align*}


\begin{thebibliography}{99}


%


\bibitem[Zou and Serfling  (2000)]{Serfling2000}
Zou, Y. and Serfling, R. General notions of statistical depth function. {\em The Annals of Statistics} {\bf 2000}, {\em 28}, 461-482.


%
%
%
%
%
%
%
%
  
  
\bibitem[Chenouri and Small  (2012)]{Small2011}
Chenouri, S. and Small, C. G. A nonparametric multivariate multisample test based on data depth. {\em Electronic Journal of Statistics} {\bf 2012}, {\em 6}, 760--782.

\bibitem[Liu and Singh  (1993)]{liu1993quality}
Liu, R. Y. and Singh, K. A quality index based on data depth and multivariate rank tests. {\em Journal of the American Statistical Association} {\bf 1993}, {\em 88(421)}, 252-260.


\bibitem[Brown  (1958)]{Brown58}
Brown,M., B. Statistical use of spatial median. {\em J.Roy.Statist.Soc.} {\bf 1958}, {\em 53}, 448--456.


\bibitem[Gower  (1974)]{Gower74}
Gower, C.,J. Algorithm as 78: The mediancentre. {\em App.Statist.} {\bf 1974}, {\em 23}, 466--470.

\bibitem[Liu  (1992)]{Liu92}
Liu, R. Y. Data depth and multivariate rank tests. {\em In $L_{1}$-Statistics and Related Methods (Y. Dodge, ed.)} {\bf 1992}, 279-294.


\bibitem[Zou  and He  (2006)]{zuo2006limiting}
Zou, Y. and He, X. One the limiting distributions of multivariate depth-based rank sum statistics and related tests. {\em The Annals of Statistics} {\bf 2006}, {\em 24(6)}, 2879--2896.

%
\bibitem[Hair et.al (1998)]{Hair}
Hair, J. F., Anderson, R. E., Tatham, R. L., and Black, W. C. Multivariate data analysis (5th ed.). {\em New York: Macmillan} {\bf 1998} Chapter 6


\bibitem[Barale, M. (2021)]{Barale&Shirke}
Barale, M.; Shirke, D. A test based on data depth for testing location-scale of the two multivariate populations. {\em Journal of statistical Computation and Simulation} {\bf 2021}, {\em 91(4)}, 768--785.


\bibitem[Liu et al. (1999)]{Liu99}
Liu, R. Y., Jesse, M. P.  and Kesar, S. Multivariate analysis by data depth: Descriptive statistics, graphics and inference. {\em The Annals of Statistics} {\bf 1999}, 783-858.

\bibitem[Li and Liu (2004)]{Liu04}
Li, J. and Liu, R.Y. New Nonparametric Tests of Multivariate Locations and Scales Using Data Depth. {\em Statistical Science} {\bf 2004}, {\em 19(4)}, 686–696.

\bibitem[Shi, Zhang and Fu (2023)]{Shi2023}
Shi, X., Zhang, Y., and Fu, Y. Two-sample tests based on data depth. {\em Entropy} {\bf 2023}, {\em 25(2)}, 238.

\bibitem[Anderson (1962)]{Anderson1962}
Anderson, T. On the Distribution of the Two-Sample Cram\'er-von Mises Criterion. {\em The Annals of Mathematical Statistics} {\bf 1962}, {\em 33(3)}, 1148–1159.

\bibitem[Kim, Balakrishnan and Wasserman (2020)]{kim2020}
Kim, I., Balakrishnan, S., and Wasserman, L. Robust multivariate nonparametric tests via projection averaging. {\em The Annals of Mathematical Statistics} {\bf 2020}, {\em 48(6)}, 3417 - 3441.
  
\bibitem[Sz\'ekely and Rizzo (2013)]{szekely2013}
Sz\'ekely , G. J. and Rizzo, M. L. Energy statistics: A class of statistics based on distances.. {\em Journal of Statistical Planning and Inference} {\bf 2013}, {\em 143(8)}, 1249–1272.

\end{thebibliography}
\end{document}